\documentclass{amsart}

\usepackage{amssymb,latexsym,amsfonts,amsmath}
\usepackage{graphicx,psfrag,epsfig,epsf}
\include{diagrams}
\usepackage{eso-pic}
\usepackage{graphicx}
\usepackage{color}
\usepackage{diagrams}
\usepackage{tikz}
\usetikzlibrary{automata}

\usepackage{amsfonts}
\usepackage{amssymb,latexsym,amsfonts,amsmath}
\usepackage{graphicx}

\newtheorem{theorem}{Theorem}[section]

\theoremstyle{definition}
\newtheorem{definition}[theorem]{Definition}
\newtheorem{example}[theorem]{Example}

\theoremstyle{remark}

\numberwithin{equation}{section}

\begin{document}

\begin{abstract}                          
Control systems are usually modeled by differential equations describing how physical phenomena can be influenced by certain control parameters or inputs. Although these models are very powerful when dealing with physical phenomena, they are less suitable to describe software and hardware interfacing the physical world. For this reason there is a growing interest in describing control systems through \emph{symbolic models} that are abstract descriptions of the continuous dynamics, where each ``symbol'' corresponds to an ``aggregate'' of states in the continuous model. Since these symbolic models are of the same nature of the models used in computer science to describe software and hardware, they provide a unified language to study problems of control in which software and hardware interact with the physical world. Furthermore the use of symbolic models enables one to leverage techniques from supervisory control and algorithms from game theory for controller synthesis purposes. In this paper we show that every incrementally globally asymptotically stable nonlinear control system is approximately equivalent (bisimilar) to a symbolic model. The approximation error is a design parameter in the construction of the symbolic model and can be rendered as small as desired. Furthermore if the state space of the control system is bounded the obtained symbolic model is finite. For digital control systems, and under the stronger assumption of incremental input--to--state stability, symbolic models can be constructed through a suitable quantization of the inputs.
\end{abstract}

\title[Approximately bisimilar symbolic models for nonlinear control systems]{Approximately bisimilar symbolic models \\
for nonlinear control systems}
\thanks{This work has been partially supported by the National Science Foundation CAREER award 0717188 and by 
the ANR SETIN project VAL-AMS.}

\author[Giordano Pola, Antoine Girard and Paulo Tabuada]{Giordano Pola$^{1}$, Antoine Girard$^{2}$ and Paulo Tabuada$^{1}$}
\address{$^{1}$Department of Electrical Engineering\\
University of California at Los Angeles,
Los Angeles, CA 90095}
\email{\{pola,tabuada\}@ee.ucla.edu}
\urladdr{http://www.ee.ucla.edu/$\sim$pola}
\urladdr{http://www.ee.ucla.edu/$\sim$tabuada}
\address{$^{2}$Universit\'e Joseph Fourier, Laboratoire Jean Kuntzmann, B.P. 53, 38041, Grenoble, France}
\email{Antoine.Girard@imag.fr}
\urladdr{http://ljk.imag.fr/membres/Antoine.Girard/}

\maketitle

\section{Introduction}
The idea of using models at different levels of abstraction has been 
 successfully used in the formal methods community with the purpose of mitigating the complexity of software verification.  A central notion when dealing with complexity reduction, is the one
of bisimulation equivalence, introduced by Milner \cite{Milner} and Park
\cite{Park} in the 80s'. The key idea is to find
and compute an equivalence relation on the state space of the system, that respects the system dynamics. This equivalence relation induces a new system on the quotient space that shares most properties of interest with the original model.
This approach leads to an alternative methodology for the analysis and control of large--scale control systems.
In fact from the analysis point of view, symbolic models provide a unified framework for describing continuous systems as well as, hardware and software interacting with the physical environment.
Furthermore, the use of symbolic models allows one to leverage the rich literature on
supervisory control \cite{RamWonI} and algorithmic approaches to game theory \cite{AVW03}, for controller design. \\
After the pioneering work of Alur and Dill \cite{TheoryTA} that showed existence of symbolic models for timed automata, researchers tried to identify more general classes of continuous systems admitting finite bisimulations. The existing results can be roughly classified into four main different lines of research:
\begin{itemize}
\item [(i)]\textbf{Simulation/bisimulation}: symbolic models have been studied in \cite{LTLControl,SMC,Girard_HSCC07} for discrete--time control systems, in \cite{Paulo_HSCC07} for continuous--time control systems and in \cite{OMinimal} for \mbox{o-minimal} hybrid systems among others. Reduction of continuous control systems to continuous control systems with lower dimensional state space has been addressed in \cite{BisimSchaft,Grasse:SIAM:07,BisimSCL,Pola:CDC04};
\item [(ii)] \textbf{Quantized control systems}: finite abstractions have been studied in \cite{BMP02,BMP06} for
certain classes of control systems with quantized inputs;
\item [(iii)] \textbf{Qualitative reasoning}: symbolic models were constructed using methods of qualitative reasoning in \cite{QR_HSCC03,QRBook};
\item [(iv)] \textbf{Stochastic automata}: abstractions of continuous--time control systems by means of stochastic automata have been studied in \cite{Lunze2001,LunzeBook}.
\end{itemize}
We defer to the last section of the paper a comparison between the results presented in this paper and the above lines of research.
In this paper we follow the line of research based on simulation/bisimulation by making use of the 
recently introduced notion of  \textit{approximate bisimulation}~\cite{AB-TAC07}, that captures equivalence of systems in an approximate setting. By relaxing the usual notion of bisimulation to approximate bisimulation, a larger class of control systems can be expected to admit symbolic models. In fact the work in \cite{Paulo_HSCC07} shows that for every asymptotically stabilizable control system it is possible to construct a symbolic model, which is based on an approximate notion of
simulation (one--sided version of bisimulation).
However, if a controller fails to exist for the symbolic model, nothing can be concluded regarding the existence of a controller for the original model. This drawback is a
direct consequence of the one--sided notion used in \cite{Paulo_HSCC07}. For this reason, an extension of the results in \cite{Paulo_HSCC07} from simulation to bisimulation is needed. The aim of this paper is precisely to provide such extension. The key idea in the results that we propose is to replace the assumption of asymptotic stabilizability of~\cite{Paulo_HSCC07} with the stronger notion of asymptotic stability. \emph{We show that every
incrementally globally asymptotically stable
nonlinear control system admits a symbolic model that is an approximate bisimulation, with a precision that is a--priori
defined, as a design parameter. Furthermore, if the state space of the control
system is bounded the symbolic
model is finite. Moreover, for incrementally input--to--state stable digital control systems, i.e. systems where control
signals are piecewise--constant, a symbolic model can
 be obtained by quantizing the space of inputs}. 
 As an illustrative example, we apply the proposed techniques to a control design problem for a pendulum.
A preliminary version of these results appeared in \cite{PolaTabuadaCDC07a}. 

\section{Control systems and stability notions\label{sec3}}
\subsection{Notations} 
The symbols $\mathbb{N}$, $\mathbb{Z}$, $\mathbb{R}$, $\mathbb{R}^{+}$ and $\mathbb{R}_{0}^{+}$ denote the natural, integers, real, positive and nonnegative real numbers, respectively. 
Given a vector $x\in\mathbb{R}^{n}$ we denote by $x'$ the transpose of $x$ and by $x_{i}$ the $i$--th element of $x$; furthermore $\Vert x\Vert$ denotes the infinity norm of $x$; we recall that 
\mbox{$\Vert x\Vert:=max\{|x_1|,|x_2|,...,|x_n|\}$}, where $|x_i|$ is the absolute value of $x_i$. 
%
%
The symbol $\mathcal{B}_{\varepsilon}(x)$ denotes the closed ball centered at
\mbox{$x\in{\mathbb{R}}^{n}$} with radius $\varepsilon\in\mathbb{R}^{+}_{0}$, i.e. \mbox{$\mathcal{B}_{\varepsilon}(x)=\{y\in{\mathbb{R}}^{n} \,:\,\Vert x-y\Vert\leq\varepsilon\}$}. For any $A\subseteq
\mathbb{R}^{n}$ and \mbox{$\mu\in{\mathbb{R}^{+}}$} define \mbox{$[A]_{\mu}:=\{a\in A\,\,|$} \mbox{$a_{i}%
=k_{i}\mu,$}\mbox{$\,\,\,k_{i}\in\mathbb{Z},$} \mbox{$i=1,...,n\}$}. The set $[A]_{\mu}$ will be used in the subsequent developments as an approximation of the set $A$ with precision $\mu$. By geometrical considerations on the infinity norm, for any $\mu\in \mathbb{R}^{+}$ and $\lambda\geq\mu/2$ the collection of sets $\{\mathcal{B}_{\lambda}(q)\}_{q\in[\mathbb{R}^{n}]_{\mu}}$ is a covering of $\mathbb{R}^{n}$, i.e.
$\mathbb{R}^{n}\subseteq{\textstyle\bigcup\nolimits_{q\in\lbrack\mathbb{R}^{n}]_{\mu}}}\mathcal{B}_{\lambda}(q)$; conversely for any $\lambda<\mu/2$, \mbox{$\mathbb{R}^{n}\nsubseteq{\textstyle\bigcup\nolimits_{q\in\lbrack\mathbb{R}^{n}]_{\mu}}}\mathcal{B}_{\lambda}(q)$}.\\ 
We now recall from \cite{Khalil,Sontag} some notions that will be employed in Sections \ref{sec:2.2} and \ref{sec:stab} to define trajectories and some stability notions for control systems.
A function $f:[a,b]\rightarrow\mathbb{R}^{n}$ is said to be absolutely continuous on $[a,b]$ if for any $\varepsilon\in \mathbb{R}^{+}$ there exists $\delta\in \mathbb{R}^{+}$ so 
that for every $k\in\mathbb{N}$ and for every sequence of points $a\leq a_{1}<a_{1}<b_{1}<a_{2}<b_{2}<...<a_{k}<b_{k}\leq b$, if $\sum_{i=1}^{m}(b_{i}-a_{i})<\delta$ then $\sum_{i=1}^{m}|f(b_{i})-f(a_{i})|<\varepsilon$. A function $f:]a,b[\rightarrow\mathbb{R}^{n}$ is said to be locally absolutely continuous if the restriction of $f$ to any compact subset of $]a,b[$ is absolutely continuous. Given a measurable function \mbox{$f:\mathbb{R}_{0}^{+}\rightarrow\mathbb{R}^{n}$}, the
\mbox{(essential)} supremum of $f$ is denoted by $\Vert f\Vert_{\infty}$; we recall that $\Vert f\Vert_{\infty}:=(ess)sup\{\Vert f(t)\Vert,$ $t\geq0\}$; 
$f$ is essentially bounded if $\Vert f \Vert_{\infty} < \infty$. 
For a given time $\tau\in\mathbb{R}^{+}$, define $f_{\tau}$ so that
$f_{\tau}(t)=f(t)$, for any $t\in [0,\tau)$, and $f(t)=0$ elsewhere; 
$f$ is said to be locally essentially bounded if for any $\tau\in\mathbb{R}^{+}$,
$f_{\tau}$ is essentially bounded.
A function \mbox{$f:\mathbb{R}^{n}\rightarrow \mathbb{R}$} is said to be radially unbounded if $f(x)\rightarrow \infty$ as $\Vert x\Vert\rightarrow \infty$. A continuous function $\gamma:\mathbb{R}_{0}^{+}%
\rightarrow\mathbb{R}_{0}^{+}$, is said to belong to class $\mathcal{K}$ if it
is strictly increasing and \mbox{$\gamma(0)=0$}; $\gamma$ is said to belong to class
$\mathcal{K}_{\infty}$ if \mbox{$\gamma\in\mathcal{K}$} and $\gamma(r)\rightarrow
\infty$ as $r\rightarrow\infty$. A continuous function \mbox{$\beta:\mathbb{R}_{0}^{+}\times\mathbb{R}_{0}^{+}\rightarrow\mathbb{R}_{0}^{+}$} is said to
belong to class $\mathcal{KL}$ if for each fixed $s$, the map $\beta(r,s)$
belongs to class $\mathcal{K}_{\infty}$ with respect to $r$ and, for each
fixed $r$, the map $\beta(r,s)$ is decreasing with respect to $s$ and
$\beta(r,s)\rightarrow0$ as \mbox{$s\rightarrow\infty$}. 
The following notions will be used in Sections \ref{sec4}, \ref{sec5} and \ref{sec6} to define the concept of approximate bisimulation and the symbolic models that we propose in this paper.
The identity map on a set $A$ is denoted by $1_{A}$. 
Given two sets $A$ and $B$, if $A$ is a subset of $B$
we denote by \mbox{$\imath_{A}:A\hookrightarrow B$} or simply by $\imath$ the natural
inclusion map taking any $a\in A$ to \mbox{$\imath (a) = a \in B$}. Given a function $f:A\rightarrow B$ the symbol $f(A)$ denotes
the image of $A$ through $f$, i.e. $f(A):=\{b\in B:\exists a\in A$ s.t.
$b=f(a)\}$. We identify a relation
$R\subseteq A\times B$ with the map \mbox{$R:A\rightarrow2^{B}$} defined by
$b\in R(a)$ if and only if \mbox{$(a,b)\in R$}. Given a relation $R\subseteq A\times B$,
$R^{-1}$ denotes the inverse relation of $R$, i.e. \mbox{$R^{-1}:=\{(b,a)\in B\times A:(a,b)\in R\}$}.\\

\subsection{Control Systems}\label{sec:2.2}

The class of control systems that we consider in this paper is formalized in
the following definition.

\begin{definition}
\label{Def_control_sys}A \textit{control system} is a quadruple 
\mbox{$\Sigma=(\mathbb{R}^{n},U,\mathcal{U},f)$},
where:
\begin{itemize}
\item $\mathbb{R}^{n}$ is the state space;
\item $U\subseteq\mathbb{R}^{m}$ is the input space;
\item $\mathcal{U}$ is a subset of the set of all locally essentially bounded functions of time from
intervals of the form $]a,b[\subseteq\mathbb{R}$ to $U$ with $a<0$ and $b>0$;
\item $f:\mathbb{R}^{n}\times U\rightarrow\mathbb{R}^{n}$ is a continuous map
satisfying the following Lipschitz assumption: for every compact set
$K\subset\mathbb{R}^{n}$, there exists a constant $\kappa>0$ such that $\Vert f(x,u)-f(y,u)\Vert\leq \kappa\Vert x-y\Vert$, for all $x,y\in K$ and all $u\in U$.
\end{itemize}
\end{definition}
A locally absolutely continuous curve $\mathbf{x}:]a,b[\rightarrow\mathbb{R}^{n}$ is said to be a
\textit{trajectory} of $\Sigma$ if there exists $\mathbf{u}\in\mathcal{U}$
satisfying 
$\dot{\mathbf{x}}(t)=f(\mathbf{x}(t),\mathbf{u}(t))$, 
for almost all $t\in$ $]a,b[$. 
Although we have defined trajectories over open domains, we shall refer to
trajectories $\mathbf{x:}[0,\tau]\rightarrow\mathbb{R}^{n}$ defined on closed
domains $[0,\tau],$ $\tau\in\mathbb{R}^{+}$ with the understanding of the
existence of a trajectory $\mathbf{z}:]a,b[\rightarrow\mathbb{R}^{n}$
such that $\mathbf{x}=\mathbf{z}|_{[0,\tau]}$. We will also write
$\mathbf{x}(t,x,\mathbf{u})$ to denote the point reached at time $t\in]a,b[$
under the input $\mathbf{u}$ from initial condition $x$; this point is
uniquely determined, since the assumptions on $f$ ensure existence and
uniqueness of trajectories \cite{Sontag}.\\
A control system $\Sigma$ is said to be
\textit{forward complete} if every trajectory is defined on an interval of the
form $]a,\infty\lbrack$. Sufficient and necessary conditions for a system to be forward complete can be found in \cite{fc-theorem}.
Simpler, but only sufficient, conditions for forward completeness are also available in the literature. These include linear growth or compact support of the vector field (see e.g. \cite{OCT}).
\subsection{Stability notions}\label{sec:stab}
The results presented in this paper will assume certain stability assumptions that we briefly recall in this section.
\begin{definition}
\label{Def_IRAS} \cite{IncrementalS} 
A control system $\Sigma$ is 
\textit{incrementally globally asymptotically stable (}$\delta$--GAS) if it is
forward complete and there exist a $\mathcal{KL}$ function $\beta$ such that
for any $t\in\mathbb{R}^{+}_{0}$, any $x,y\in\mathbb{R}^{n}$ and any $\mathbf{u}\in\mathcal{U}$ the following condition is satisfied:%
\begin{equation}
\left\Vert \mathbf{x}(t,x,\mathbf{u})-\mathbf{x}(t,y,\mathbf{u}%
)\right\Vert \leq\beta(\left\Vert x-y\right\Vert ,t).
\label{deltaUGAS}%
\end{equation}
\end{definition}
Definition above can be thought of as an incremental version of the classical notion of
global asymptotic stability (GAS) \cite{Khalil}. 
\begin{definition}
\label{Def_IISS} \cite{IncrementalS} 
A control system $\Sigma$ is 
\textit{incrementally input--to--state stable} ($\delta$--ISS) if it is
forward complete and there exist a $\mathcal{KL}$ function $\beta$ and a
$\mathcal{K}_{\infty}$ function $\gamma$ such that for any $t\in\mathbb{R}^{+}_{0}$, any
$x,y\in\mathbb{R}^{n}$ and any $\mathbf{u}%
,\mathbf{v}\in\mathcal{U}$ the following condition is satisfied:
\begin{eqnarray}
\left\Vert \mathbf{x}(t,x,\mathbf{u})-\mathbf{x}(t,y%
,\mathbf{v})\right\Vert & \leq &\beta(\left\Vert x-y\right\Vert
,t)+\gamma(\left\Vert \mathbf{u}-\mathbf{v}\right\Vert _{\infty}).\notag \\
\label{deltaISS}
\end{eqnarray}
\end{definition}
It is readily seen, by observing
(\ref{deltaUGAS}) and (\ref{deltaISS}), 
that
$\delta$--ISS implies $\delta$--GAS, while the converse is not true in
general (see \cite{IncrementalS} for some examples). \\
In general, inequalities (\ref{deltaUGAS}) and (\ref{deltaISS}) are difficult
to check directly. Fortunately $\delta$--GAS and $\delta$--ISS can be characterized by dissipation inequalities. %
\begin{definition}
\label{LyapunovGAS}
Consider a control system $\Sigma$ and a smooth function
\mbox{$V:\mathbb{R}^{n}\times\mathbb{R}^{n}\rightarrow\mathbb{R}_{0}^{+}$}. Function $V$ is called a $\delta$--GAS\textit{ Lyapunov function} for $\Sigma$, if there exist
$\mathcal{K}_{\infty}$ functions $\alpha_{1}$, $\alpha_{2}$ and $\rho$\ such that:
\\
\\
(i)  for any $x,y\in\mathbb{R}^{n}$%
\[
\alpha_{1}(\left\Vert x-y\right\Vert )\leq V(x,y)\leq
\alpha_{2}(\left\Vert x-y\right\Vert )\text{;}
\]
(ii) for any $x,y\in\mathbb{R}^{n}$ and any $u\in U$%
\[
\frac{\partial V}{\partial x}f(x,u)+\frac{\partial V}{\partial y%
}f(y,u)<-\rho(\left\Vert x-y\right\Vert )\text{.}%
\] 
Function $V$ is
called a $\delta$\textit{--ISS Lyapunov function} for $\Sigma$, if there exist
$\mathcal{K}_{\infty}$ functions $\alpha_{1}$, $\alpha_{2}$, $\rho$ and
$\sigma$\ satisfying conditions (i) and:
\\
\\
(iii) for any $x,y\in\mathbb{R}^{n}$ and any $u,v\in U$
\[
\frac{\partial V}{\partial x}f(x,u)+\frac{\partial V}{\partial
y}f(y,v) <-\rho(\left\Vert x-y\right\Vert )+\sigma
(\left\Vert u-v\right\Vert )\text{.}
\]
\end{definition}
The following result completely characterizes $\delta$--GAS and $\delta$--ISS
in terms of existence of Lyapunov functions.
\begin{theorem}
\cite{IncrementalS}
Consider a control system $\Sigma=(\mathbb{R}^{n},U,\mathcal{U},f)$. Then:
\begin{itemize}
	\item If $U$ is compact then $\Sigma$ is $\delta$--GAS if and only if it admits a $\delta$--GAS Lyapunov function;
	\item If $U$ is closed, convex, contains the origin and \mbox{$f(0,0)=0$}, then $\Sigma$ is $\delta$--ISS
 if it admits a $\delta$--ISS Lyapunov
function. Moreover if $U$ is compact, existence of a $\delta$--ISS Lyapunov
function is equivalent to $\delta$--ISS.
\end{itemize}
\label{TH-IISS}
\end{theorem}
%
\section{Approximate bisimulation\label{sec4}}
In this section we introduce a notion of approximate equivalence upon which
all the results in this paper rely. We start by introducing the class of transition systems that will be used in
this paper as abstract models for control systems.
\begin{definition}
A transition system is a quintuple 
\mbox{$T=(Q,L,\rTo,O,H)$}, 
consisting of:
\begin{itemize}
\item A set of states $Q$;
\item A set of labels $L$;
\item A transition relation $\rTo\subseteq Q\times L\times Q$;
\item An output set $O$;
\item An output function $H:Q\rightarrow O$.
\end{itemize}
A transition system $T$ is said to be: 
\begin{itemize}
	\item \textit{metric}, if the output set $O$ is equipped with a metric \mbox{$\mathbf{d}:O\times O\rightarrow\mathbb{R}_{0}^{+}$};
	\item \textit{countable}, if $Q$ and $L$ are countable sets;
	\item \textit{finite}, if $Q$ and $L$ are finite sets. 
\end{itemize}
\end{definition}
We will follow standard practice and denote an element $(q,l,p)\in
\rTo$ by $q\rTo^{l} p$. Transition systems
capture dynamics through the transition relation. For any states \mbox{$q,p\in Q$},
$q\rTo^{l} p$ simply means that it is possible to evolve or
jump from state $q$ to state $p$ under the action labeled by $l$. 
We will use transition systems as an abstract representation of control
systems. There are several different ways in which control
systems can be transformed into transition systems. We now describe one of these, which has the
property of capturing all the information contained in a control system
$\Sigma$. \\
Given a control system $\Sigma=(\mathbb{R}^{n},U,\mathcal{U},f)$ define the
transition system:
\begin{equation}
T(\Sigma):=(Q,L,\rTo,O,H),
\label{T0}
\end{equation}
where:
\begin{itemize}
\item $Q=\mathbb{R}^{n}$;
\item $L=\mathcal{U}$;
\item $q\rTo^{\mathbf{u}} p$, if 
$\mathbf{x}(\tau,q,\mathbf{u})=p$ for some $\tau\in\mathbb{R}^{+}$;
\item $O=\mathbb{R}^{n}$;
\item $H=1_{\mathbb{R}^{n}}$.
\end{itemize}
Transition system $T(\Sigma)$ is metric when we regard the set
\mbox{$O=\mathbb{R}^{n}$} as being equipped with the metric \mbox{$\mathbf{d}(p,q)=\Vert p-q\Vert$}. Note that the state space of $T(\Sigma)$ is infinite. The aim of this paper is to study existence of countable transition systems that are approximately equivalent to $T(\Sigma)$. The notion of 
equivalence that we consider is the one of \textit{bisimulation equivalence}
\cite{Milner,Park}. Bisimulation relations are standard mechanisms to relate
the properties of transition systems
. Intuitively, a
bisimulation relation between a pair of transition systems $T_{1}$ and $T_{2}$
is a relation between the corresponding state sets explaining how a state trajectory $r_{1}$ of $T_{1}$ can be transformed into a state 
trajectory $r_{2}$ of $T_{2}$ and vice versa. While typical bisimulation
relations require that $r_{1}$ and $r_{2}$ are observationally
indistinguishable, that is $H_{1}(r_{1})=H_{2}(r_{2})$, we shall relax this
by requiring $H_{1}(r_{1})$ to simply be close to $H_{2}(r_{2})$ where
closeness is measured with respect to the metric on the output set. The
following notion has been introduced in \cite{AB-TAC07} and in a
slightly different formulation in \cite{Paulo_HSCC07}.
\begin{definition}
\label{ASR}Let $T_{1}=(Q_{1},L_{1},\rTo_{1},O,H_{1})$ and
\mbox{$T_{2}=(Q_{2},L_{2},\rTo_{2},O,H_{2})$} be metric transition systems
with the same output set and metric $\mathbf{d}$, and let $\varepsilon \in \mathbb{R}_{0}^{+}$ 
be a given precision.
A relation $R\subseteq Q_{1}\times Q_{2}$ is said to be an $\varepsilon
$--\textit{approximate bisimulation} relation between $T_{1}$ and $T_{2}$, if for
any $(q_{1},q_{2})\in R$:\\
\\
(i)  $\mathbf{d}(H_{1}(q_{1}),H_{2}(q_{2}))\leq\varepsilon$;\\
(ii) $q_{1}\rTo^{l_{1}}_{1} p_{1}$ implies 
existence of $q_{2}\rTo^{l_{2}}_{2} p_{2}$ such that
$(p_{1},p_{2})\in R$.\\
(iii) $q_{2}\rTo^{l_{2}}_{2} p_{2}$ implies 
existence of $q_{1}\rTo^{l_{1}}_{1} p_{1}$ such that
$(p_{1},p_{2})\in R$.
\\
Moreover $T_{1}$ is $\varepsilon$\textit{--bisimilar} to $T_{2}$ if there exists
an \mbox{$\varepsilon$--approximate} bisimulation relation $R$\ between $T_{1}$ and $T_{2}$ such that
\mbox{$R(Q_{1})=Q_{2}$} and $R^{-1}(Q_{2})=Q_{1}$. 
\end{definition}
%
\section{Approximate bisimilar symbolic models\label{sec5}}
In the following we will work with a sub--transition system of $T(\Sigma)$
obtained by selecting those transitions from $T(\Sigma)$ that describe
trajectories of duration $\tau$ for some chosen \mbox{$\tau\in\mathbb{R}^{+}$}. This
can be seen as a time discretization or sampling process. Given a control system $\Sigma$ and a
parameter \mbox{$\tau\in\mathbb{R}^{+}$} define the transition system:%
\[
T_{\tau}(\Sigma):=(Q_{1},L_{1},\rTo_{1},O_{1},H_{1}),
\]
where:
\begin{itemize}
\item $Q_{1}=\mathbb{R}^{n}$;
\item $L_{1}=\{l_{1}\in\mathcal{U}\,\,\vert\,\,\mathbf{x}(\tau
,x,l_{1})$ is defined for all $x\in\mathbb{R}^{n}\}$;
\item $q\overset{l_{1}}{\rTo_{1}}p$, if 
$\mathbf{x}(\tau,q,l_{1})=p$;
\item $O_{1}=\mathbb{R}^{n}$;
\item $H_{1}=1_{\mathbb{R}^{n}}$.
\end{itemize}
Transition system $T_{\tau}(\Sigma)$ is metric when we regard
\mbox{$O_{1}=\mathbb{R}^{n}$} as being equipped with the metric \mbox{$\mathbf{d}(p,q)=\Vert p-q\Vert$}.
Note that the set of labels $L_{1}$ is composed by (only) those control signals of $\mathcal{U}$ for which a trajectory of $\Sigma$ exists for any time $t\in[0,\tau]$
and \textit{for any initial condition}
$x\in\mathbb{R}^{n}$. 
Any measurable control input can be included in $L_{1}$ when the control system is forward complete.\\
In the following we show existence of a countable transition system that is
 approximately bisimilar to $T_{\tau}(\Sigma)$, provided that $\Sigma$
satisfies some stability properties.\\
By simple considerations on the infinity norm, for any given precision $\eta\in\mathbb{R}^{+}$ we can
approximate the state space $Q_{1}=\mathbb{R}^{n}$ of $T_{\tau}(\Sigma)$ by means of the
countable set $Q_{2}:=[\mathbb{R}^{n}]_{\eta}$ so that for any $x\in\mathbb{R}^{n}$ 
there exists $q\in Q_{2}$ such that $\Vert x-q \Vert \leq \eta/2$.\\
The approximation of the set of labels $L_{1}$ of $T_{\tau}(\Sigma)$ is more involved.
We approximate $L_{1}$ by means of the set:
\begin{equation}
L_{2}:=%
{\textstyle\bigcup\nolimits_{q\in Q_{2}}}
L_{2}(q),
\label{L2}
\end{equation}
where $L_{2}(q)$ captures the set of labels that can be applied at the state $q\in Q_{2}$ of the symbolic model. The definition of $L_{2}(q)$ is based on the notion of reachable sets. 
Given any 
state $q\in Q_{1}$ consider the set:
\begin{equation}
\mathcal{R}(\tau,q)=\left\{p\in Q_{1}:q\rTo^{l_{1}}_{1}p,l_{1}\in L_{1}\right\},
\label{reach}
\end{equation}
of reachable states of $T_{\tau}(\Sigma)$ from $q$. Notice that 
$\mathcal{R}(\tau,q)$ is well defined because of the definition of the set of labels $L_{1}$.
We approximate $\mathcal{R}(\tau,q)$ by means of a countable set, as follows. Given any precision $\mu\in\mathbb{R}^{+}$, consider the set:%
\[
\mathcal{P}_{\mu}(\tau,q):=\{y\in\lbrack\mathbb{R}^{n}]_{\mu}:\exists
z\in\mathcal{R}(\tau,q)\text{ s.t. }\Vert y-z\Vert\leq\mu/2\}\text{,}%
\]
and define the function 
$\psi_{\mu}^{\tau,q}:\mathcal{P}_{\mu}(\tau,q)\rightarrow L_{1}$, %
that associates to any $y\in\mathcal{P}_{\mu}(\tau,q)$ a label
\mbox{$l_{1}=\psi_{\mu}^{\tau,q}(y)\in L_{1}$} so
that $\Vert y-\mathbf{x}(\tau,q,l_{1})\Vert\leq\mu/2$. Notice that the function
$\psi_{\mu}^{\tau,q}$ is not unique. The set $L_{2}(q)$ appearing in (\ref{L2}) can now be defined by 
$L_{2}(q):=\psi_{\mu}^{\tau,q}(\mathcal{P}_{\mu}(\tau,q))$.
Notice that since $L_{2}(q)$ is the image through $\psi_{\mu}^{\tau,q}$ of a countable set, it is countable. Therefore $L_{2}$ as defined in (\ref{L2}) is countable, as well. Furthermore the set $L_{2}$ approximates the set $L_{1}$ in the sense that given any $q\in Q_{2}$, for any $l_{1}\in L_{1}$ there exists $l_{2}\in L_{2}(q)$ so that: 
\begin{equation}
\Vert \mathbf{x}(\tau,q,l_{1})-\mathbf{x}(\tau,q,l_{2})\Vert\leq\mu.
\label{label}
\end{equation}
%
%
We now have all the ingredients to define a symbolic model that will be used to approximate a control system.\\
Given a control system $\Sigma=(\mathbb{R}^{n},U,\mathcal{U},f)$, any $\tau\in\mathbb{R}^{+}$, $\eta
\in\mathbb{R}^{+}$ and $\mu\in\mathbb{R}^{+}$\ define the following transition
system:
\begin{equation}
T_{\tau,\eta,\mu}(\Sigma):=(Q_{2},L_{2},\rTo_{2},O_{2},H_{2}),
\label{T1}%
\end{equation}
where:
\begin{itemize}
\item $Q_{2}=[\mathbb{R}^{n}]_{\eta};$
\item $L_{2}=%
{\textstyle\bigcup\nolimits_{q\in Q_{2}}}
L_{2}(q)$;
\item $q\rTo^{l}_{2} p$, if $l\in
L_{2}(q)$ and $\left\Vert p-\mathbf{x}(\tau,q,l)\right\Vert \leq\eta/2$;
\item $O_{2}=\mathbb{R}^{n}$;
\item $H_{2}=\imath : Q_{2} \hookrightarrow O_{2}$.
\end{itemize}
We think of $T_{\tau,\eta,\mu}(\Sigma)$ as a metric transition system where $O_{2}=\mathbb{R}^{n}$ is equipped with the metric \mbox{$\mathbf{d}(p,q)=\Vert p-q\Vert$}. Parameters $\tau\in\mathbb{R}^{+}$, $\eta\in\mathbb{R}^{+}$ and $\mu
\in\mathbb{R}^{+}$ in $T_{\tau,\eta,\mu}(\Sigma)$ can be
thought of, respectively, as a sampling time, a state space and an input space
quantization. \\
We emphasize that transition system
$T_{\tau,\eta,\mu}(\Sigma)$\ is countable because the sets $Q_{2}$
and $L_{2}$ are countable. Furthermore if the state space
of the control system $\Sigma$ is bounded, the corresponding transition system
$T_{\tau,\eta,\mu}(\Sigma)$ is \emph{finite}.\\
Note that in the definition of the transition relation $\rTo_{2}$ we require $\mathbf{x}(\tau,q,l)$ to be in the closed ball $\mathcal{B}_{\eta/2}(p)$. We can instead, require $\mathbf{x}(\tau,q,l)$ to be in $\mathcal{B}_{\lambda}(p)$ for any $\lambda\geq\eta/2$. However, we chose $\lambda=\eta/2$ because $\eta/2$ is the smallest value of $\lambda\in\mathbb{R}^{+}$ that ensures $\mathbb{R}^{n}\subseteq{\textstyle\bigcup\nolimits_{q\in\lbrack\mathbb{R}^{n}]_{\eta}}}\mathcal{B}_{\lambda}(q)$. In fact, this choice of $\lambda$ reduces the number of transitions in the definition of the symbolic model in (\ref{T1}).\\
%
We can now give the main result of this paper which relates $\delta
$--GAS to existence of symbolic model.
\begin{theorem}
\label{Th_main}Consider a control system $\Sigma$ and any desired precision $\varepsilon\in\mathbb{R}^{+}$.
If $\Sigma$ is $\delta$--GAS then for any $\tau\in\mathbb{R}^{+}$, $\eta
\in\mathbb{R}^{+}$ and $\mu\in\mathbb{R}^{+}$ satisfying the following
inequality:%
\begin{equation}
\beta(\varepsilon,\tau)+\mu+\eta/2\leq\varepsilon, \label{cond}%
\end{equation}
the transition system $T_{\tau}(\Sigma)$ is $\varepsilon$--bisimilar to
$T_{\tau,\eta,\mu}(\Sigma)$.
\end{theorem}
Before giving the proof of this result we point out that if $\Sigma$ is
$\delta$--GAS, there always exist parameters $\tau\in\mathbb{R}^{+}$,
$\eta\in\mathbb{R}^{+}$ and $\mu\in\mathbb{R}^{+}$ satisfying condition
(\ref{cond}). Indeed since $\beta$ is a $\mathcal{KL}$ function, there exists a
sufficiently large value of $\tau$ so that $\beta(\varepsilon
,\tau)<\varepsilon$; then by choosing sufficiently small values of $\mu$ and
$\eta$, condition (\ref{cond}) is fulfilled.
\\
\\
\begin{proof}
Consider the relation $R\subseteq Q_{1}\times Q_{2}$ defined by $(x,q)\in R$
if and only if $||x-q||\leq\varepsilon$. By construction $R(Q_{1})=Q_{2}$; furthermore \mbox{$
Q_{1}\subseteq {\textstyle\bigcup\nolimits_{q_{2}\in Q_2}}
\mathcal{B}_{\eta/2}(q_{2})$} and therefore since by
(\ref{cond}), $\eta/2<\varepsilon$, we have that \mbox{$R^{-1}(Q_{2})=Q_{1}$}. We now
show that $R$ is an \mbox{$\varepsilon$--approximate} bisimulation relation between
$T_{\tau}(\Sigma)$ and $T_{\tau,\eta,\mu}(\Sigma)$. Consider any $(x,q)\in R$. Condition (i) in Definition \ref{ASR} is satisfied
by definition of $R$. Let us now show that condition (ii) in Definition
\ref{ASR} holds. Consider any $l_{1}\in L_{1}$ and the transition 
$x \rTo^{l_{1}}_{1} y$ in $T_{\tau}(\Sigma)$. Let $v=\mathbf{x}(\tau,q,l_{1})$; since $\mathbb{R}^{n}\subseteq {\textstyle\bigcup\nolimits_{{w}\in \lbrack\mathbb{R}^{n}]_{\mu}}}\mathcal{B}_{\mu/2}(w)$, there exists $w\in \lbrack\mathbb{R}^{n}]_{\mu}$ such that:
\begin{equation}
\Vert v-w\Vert\leq\mu/2.\label{b01}%
\end{equation}
Since $v\in {\mathcal R}(\tau,q)$, it is clear that $w\in {\mathcal P}_\mu(\tau,q)$ by definition of
${\mathcal P}_\mu(\tau,q)$.
Then, let $l_{2}\in L_2(q)$ be given by $l_{2} = \psi_{\mu}^{\tau,q}(w)$. By definition of $\psi_{\mu}^{\tau,q}$ 
and by setting $z={\bf x} (\tau,q,l_{2})$,
it follows that:
\begin{equation}
\Vert w-z\Vert\leq\mu/2.\label{b02}%
\end{equation}
Since $Q_{1}\subseteq {\textstyle\bigcup\nolimits_{q_{2}\in Q_2}}
\mathcal{B}_{\eta/2}(q_{2})$, there exists $p\in Q_2$ such that: 
\begin{equation}
\Vert z-p\Vert\leq\eta/2.\label{b2}%
\end{equation}
Thus, $q \rTo^{l_{2}}_{2} p$ in $T_{\tau,\eta,\mu}(\Sigma)$ and since $\Sigma$ is $\delta$--GAS and by (\ref{b01}), (\ref{b02}), (\ref{b2})
and (\ref{cond}), the following chain of inequalities holds:%
\begin{align*}
\Vert y-p\Vert &  =\Vert y-v+v-w+w-z+z-p\Vert\\
&  \leq\Vert y-v\Vert+\Vert v-w\Vert+\Vert w-z\Vert+\Vert z-p\Vert\\
&  \leq\beta(||x-q||,\tau)+\mu/2+\mu/2+\eta/2\\
&  \leq\beta(\varepsilon,\tau)+\mu+\eta/2\leq\varepsilon.
\end{align*}
Hence $(y,p)\in R$ and condition (ii) in Definition \ref{ASR} holds. 
We now show that also condition (iii) holds.
Consider any $(x,q)\in R$, any $l_{2}\in L_{2}$ and the transition 
$q \rTo^{l_{2}}_{2} p$ in $T_{\tau,\eta,\mu}(\Sigma)$. By
definition of $T_{\tau,\eta,\mu}(\Sigma)$: 
\begin{equation}
\left\Vert z-p\right\Vert \leq\eta/2,
\label{gio}
\end{equation}
where $z=\mathbf{x}(\tau,q,l_{2})\in Q_{1}$. Choose $l_{1}=l_{2}\in L_{1}$ and consider the transition
$x \rTo^{l_{1}}_{1} y$ in $T_{\tau}(\Sigma)$. Since
$\Sigma$ is \mbox{$\delta$--GAS} and by conditions (\ref{gio}) and (\ref{cond}), the following chain
of inequalities holds:%
\begin{align*}
\Vert y-p\Vert &  =\Vert y-z+z-p\Vert\leq\Vert y-z\Vert+\Vert z-p\Vert\\
&  \leq\beta(\Vert x-q\Vert,\tau)+\eta/2\leq\beta(\varepsilon,\tau)+\eta/2\leq\varepsilon.
\end{align*}
Thus $(y,p)\in R$, which completes the proof.
\end{proof}
Conditions of Theorem \ref{Th_main} require the control system $\Sigma$ to be \textit{globally} $\delta$--GAS as in Definition \ref{Def_IRAS}. However, it is easy to see from the above proof that this stability property can be relaxed to hold  \textit{locally}, i.e. for initial states $x,y\in\mathbb{R}^{n}$ satisfying $\Vert x-y \Vert\leq \varepsilon$.
Moreover, this stability condition is not far from also being necessary. 
The following counterexample shows that unstable control systems do not admit, in general, countable symbolic models.
\begin{example}
Consider a control system $\Sigma=(\mathbb{R},U,\mathcal{U},f)$, where $U=\{0\}$, $\mathcal{U}=\{\mathbf{0}\}$, $\mathbf{0}$ is the identically null input and $f(x)=x$. System $\Sigma$ is unstable and hence not \mbox{$\delta$--GAS}. We now show that for any $\varepsilon\in \mathbb{R}^{+}_{0}$, any $\tau\in \mathbb{R}^{+}$ and any countable transition system $T$, transition systems $T_{\tau}(\Sigma)$ and $T$ are \textit{not} $\varepsilon$--bisimilar. 
Consider any countable metric transition system $T=(Q,L,\rTo,\mathbb{R},H)$, with $H:Q\rightarrow \mathbb{R}$ and the same metric $\mathbf{d}(p,q)=\Vert p-q \Vert$ of $T_{\tau}(\Sigma)$. Consider any relation $R\subseteq Q_{1}\times Q$ satisfying conditions (i), (ii) and (iii) of Definition \ref{ASR} and such that $R(Q_{1})=Q$ and $R^{-1}(Q)=Q_{1}$. We now show that such relation $R$ does not exist.
By countability of $T$, there exist $q_{0}\in Q$ and $x_{0},y_{0}\in Q_{1}=\mathbb{R}$ such that $x_{0}\neq y_{0}$, and $(x_{0},q_{0}),(y_{0},q_{0})\in R$. Set $x_{k}=e^{\tau k}x_{0}$, $y_{k}=e^{\tau k}y_{0}$, for any $k\in\mathbb{N}$. 
Since $x_{0}\neq y_{0}$, by selecting $\lambda\in\mathbb{R}^{+}$ such that $\Vert x_{0}-y_{0} \Vert>\lambda$, we have:
\begin{equation}
\Vert x_{k}-y_{k} \Vert = e^{\tau k}\Vert x_{0}-y_{0} \Vert > e^{\tau k}\lambda,\forall k\in\mathbb{N}.
\label{zexample}
\end{equation}
Choose $k' \in \mathbb{N}$ so that $e^{\tau k'} \lambda - \varepsilon>\varepsilon$. By condition (iii) in Definition \ref{ASR} and since $R(Q_{1})=Q$ and $R^{-1}(Q)=Q_{1}$, there must exist $q_{k'}\in Q$ so that, $(x_{k'},q_{k'}),(y_{k'},q_{k'})\in R$. Since $(x_{k'},q_{k'})\in R$,
\begin{equation}
\Vert x_{k'}-H(q_{k'}) \Vert \leq \varepsilon.
\label{xexample}
\end{equation}
By combining inequalities (\ref{zexample}) and (\ref{xexample}) and by definition of $k'$, we obtain:
\begin{eqnarray}
\Vert H(q_{k'})-y_{k'} \Vert & \geq & \Vert x_{k'} - y_{k'} \Vert -\Vert x_{k'}- H(q_{k'})\Vert \notag \\
& > & e^{\tau k'} \lambda - \varepsilon>\varepsilon.
\label{zzexample}
\end{eqnarray}
Inequality (\ref{zzexample}) shows that the pair $(y_{k'},q_{k'})\in R$ does not satisfy condition (i) of Definition \ref{ASR}. Hence, there does not exist an $\varepsilon$--approximate bisimulation relation between $T_{\tau}(\Sigma)$ and $T$ and consequently  
$T_{\tau}(\Sigma)$ and $T$ are \textit{not} \mbox{$\varepsilon$--bisimilar}.
\end{example}
Theorem~\ref{Th_main} relates $T_{\tau}(\Sigma)$ to the symbolic model in (\ref{T1}), whose construction is in general difficult, since it requires the computation of reachable sets. In the next section we show that for digital control systems a symbolic model can be obtained by quantizing the input space.
\section{Digital control systems\label{sec6}}
In this section we specialize the results of the previous section to the case
of digital control systems, i.e. control systems where control signals are
piecewise--constant. In many man made systems, input signals are often physically
implemented as piecewise--constant signals and this motivates our interest in this class of systems.\\
In the following we suppose that the input space $U$ of the considered control system $\Sigma=(\mathbb{R}^{n},U,\mathcal{U},f)$ contains the origin and that it 
is a hyper rectangle of the form  
$U:= [a_1,b_1] \times [a_2,b_2] \times ... \times [a_m,b_m]$, 
for some $a_i<b_i, i=1,2,...,m$. Furthermore we suppose that control inputs are piecewise--constant; given \mbox{$\tau\in\mathbb{R}^{+}$}, 
the class of inputs that we consider is:%
\[
\mathcal{U}_{\tau}:=\{\mathbf{u}\in\mathcal{U}:\mathbf{u}(t)=\mathbf{u}%
(0),t\in\lbrack0,\tau]\}.
\]
For notational simplicity, we denote by $u$ the control input $\mathbf{u}\in
\mathcal{U}_{\tau}$ for which $\mathbf{u}(t)=u,$ $t\in\lbrack0,\tau]$. \\
Let us denote by $T_{\mathcal{U}_{\tau}}(\Sigma)$ the sub--transition system
of $T_{\tau}(\Sigma)$ where only control inputs in $\mathcal{U}_{\tau}$ are
considered. More formally define:%
\[
T_{\mathcal{U}_{\tau}}(\Sigma):=(Q_{1},L_{1},\rTo_{1},O_{1}%
,H_{1}),
\]
where:
\begin{itemize}
\item $Q_{1}=\mathbb{R}^{n}$;
\item $L_{1}=\{l_{1}\in U\,\,\vert\,\,\mathbf{x}(\tau
,x,l_{1})$ is defined for all $x\in\mathbb{R}^{n}\}$;
\item $q \rTo^{l}_{1} p$, if 
$\mathbf{x}(\tau,q,l)=p$;
\item $O_{1}=\mathbb{R}^{n}$;
\item $H_{1}=1_{\mathbb{R}^{n}}$.
\end{itemize}
Transition system $T_{\mathcal{U}_{\tau}}(\Sigma)$ is metric when we regard
$O=\mathbb{R}^{n}$ as being equipped with the metric \mbox{$\mathbf{d}(p,q)=\Vert p-q\Vert$}.
Note that analogously to $T_{\tau}(\Sigma)$, transition system $T_{\mathcal{U}_{\tau}}(\Sigma)$ is not countable. Therefore we now define a suitable countable transition system that will approximate $T_{\mathcal{U}_{\tau}}(\Sigma)$ with any desired precision.\\
Given a control system $\Sigma$, any $\tau\in\mathbb{R}^{+}$, $\eta
\in\mathbb{R}^{+}$ and $\mu\in\mathbb{R}^{+}$,\ define the following
transition system:
\begin{equation}
T_{\tau,\eta,\mu}(\Sigma):=(Q_{2},L_{2},\rTo_{2},O_{2},H_{2}),
\label{T2}%
\end{equation}
where:
\begin{itemize}
\item $Q_{2}=[\mathbb{R}^{n}]_{\eta}$;
\item $L_{2}=[L_{1}]_{\mu}$;
\item $q \rTo^{l}_{2} p$, if $\left\Vert p-\mathbf{x}(\tau,q,l)\right\Vert \leq\eta/2$;
\item $O_{2}=\mathbb{R}^{n}$;
\item $H_{2}=\imath : Q_{2} \hookrightarrow O_{2}$.
\end{itemize}
Analogously to transition system in (\ref{T1}), transition system in (\ref{T2}) is countable. Notice that transition system in (\ref{T2})
differs from the one in (\ref{T1}), (only) in the way that control inputs are approximated. In
particular, the choice of labels in transition system in (\ref{T2}) does not require
the knowledge of reachable set associated with $\Sigma$. This feature is essential when constructing the symbolic model. The computation of $\mathbf{x}(\tau,q,l)$ can be done either analytically or numerically; in the later case, numerical errors can be incorporated in the model, as follows.
Suppose there exists a parameter $\nu\in \mathbb{R}^{+}_{0}$ so that for any state $q\in Q_{2}$ and control input $l\in L_{2}$, it is possible to evaluate $\mathbf{x}(\tau,q,l)$ by means of the numerical solution $\mathbf{\tilde{x}}(\tau,q,l)$ with precision $\nu$, i.e. 
$\Vert \mathbf{x}(\tau,q,l)- \mathbf{\tilde{x}}(\tau,q,l) \Vert \leq \nu$. 
Then, the transition relation $\rTo_{2}$ in the transition system of (\ref{T2}), can be adapted to this case by requiring that $q \rTo^{l}_{2} p$, if $\left\Vert p-\mathbf{\tilde{x}}(\tau,q,l)\right\Vert\leq \eta/2-\nu$. In fact:
\begin{eqnarray}
\left\Vert p-\mathbf{x}(\tau,q,l)\right\Vert 
&
\leq 
&
\left\Vert p-\mathbf{\tilde{x}}(\tau,q,l)\right\Vert + \left\Vert \mathbf{\tilde{x}}(\tau,q,l)-\mathbf{x}(\tau,q,l)\right\Vert \notag \\
& \leq
&
\eta/2-\nu+\nu=\eta/2,\notag
\end{eqnarray}
and therefore we can recover transition relation $\rTo_{2}$, as defined in transition system (\ref{T2}).\\
We can now state the following result that relates $\delta$--ISS to the
existence of symbolic models for digital control systems.
\begin{theorem}
\label{Th:main5}Consider a control system $\Sigma$ and any desired precision $\varepsilon\in\mathbb{R}^{+}$.
If $\Sigma$ is $\delta$--ISS then for any $\tau\in\mathbb{R}^{+}$, $\eta
\in\mathbb{R}^{+}$, and $\mu\in\mathbb{R}^{+}$ satisfying the following
inequality:%
\begin{equation}
\beta(\varepsilon,\tau)+\gamma(\mu)+\eta/2\leq\varepsilon, \label{cond2}%
\end{equation}
the transition system $T_{\mathcal{U}_{\tau}}(\Sigma)$ is $\varepsilon
$--bisimilar to $T_{\tau,\eta,\mu}(\Sigma)$.
\end{theorem}
Before giving the proof of this result we point out that, analogously to
condition (\ref{cond}) of Theorem \ref{Th_main}, there always exist parameters $\tau
\in\mathbb{R}^{+}$, $\eta\in\mathbb{R}^{+}$, and $\mu\in\mathbb{R}^{+}$
satisfying condition (\ref{cond2}).
\\
\\
\begin{proof}
Consider the relation $R\subseteq Q_{1}\times Q_{2}$ defined by $(x,q)\in R$
if and only if $||x-q||\leq\varepsilon$. By construction $R(Q_{1})=Q_{2}$; since $Q_{1}\subseteq
{\textstyle\bigcup\nolimits_{q_{2}\in Q_{2}}}
\mathcal{B}_{\eta/2}(q_{2})$ and by
(\ref{cond2}), $\eta/2<\varepsilon$, we have that $R^{-1}(Q_{2})=Q_{1}$. We now
show that $R$ is an \mbox{$\varepsilon$--approximate} bisimulation relation between
$T_{\mathcal{U}_{\tau}}(\Sigma)$ and $T_{\tau,\eta,\mu}(\Sigma)$.
Consider any $(x,q)\in R$. Condition (i) in Definition \ref{ASR} is satisfied
by the definition of $R$. Let us now show that condition (ii) in Definition
\ref{ASR} holds. Consider any $l_{1}\in L_{1}$ and the transition 
$x \rTo^{l_{1}}_{1} y$ in $T_{\mathcal{U}_{\tau}}(\Sigma)$. Consider a label $l_{2}\in L_{2}$ such that:
\begin{equation}
\left\Vert l_{1}-l_{2}\right\Vert \leq\mu, \label{a5}%
\end{equation}
and set $z=\mathbf{x}(\tau,q,l_{2})$. (Notice that such label $l_{2}\in L_{2}$
exists because the assumptions on $U$ make $L_{2}=[L_{1}]_{\mu}$ non--empty.) 
For later use notice that since $l_{1}$ and $l_{2}$ are constant functions, then $\left\Vert l_{1}-l_{2}\right\Vert=\left\Vert l_{1}-l_{2}\right\Vert_{\infty}$.
Since $Q_{1}\subseteq%
{\textstyle\bigcup\nolimits_{q_{2}\in\lbrack\mathbb{R}^{n}]_{\eta}}}
\mathcal{B}_{\eta/2}(q_{2})$, there exists $p\in Q_{2}$
such that:
\begin{equation}
\left\Vert z-p\right\Vert \leq\eta/2, \label{a2}%
\end{equation}
and therefore $q \rTo^{l_{2}}_{2} p$ in $T_{\tau,\eta,\mu}(\Sigma)$. Since $\Sigma$ is $\delta$--ISS
and by (\ref{a5}), (\ref{a2}) and (\ref{cond2}), the following chain of
inequalities holds:%
\begin{eqnarray}
\Vert y-p\Vert &  =\Vert y-z+z-p\Vert\leq\Vert y-z\Vert+\Vert z-p\Vert\notag\\
&  \leq\beta(\Vert x-q\Vert,\tau)+\gamma(\Vert l_{1}-l_{2}\Vert_{\infty}%
)+\eta/2\notag\\
&  \leq\beta(\varepsilon,\tau)+\gamma(\mu)+\eta/2\leq\varepsilon.
\label{gio21}
\end{eqnarray}
Hence $(y,p)\in R$ and condition (ii) in Definition \ref{ASR} holds. 
We now show that also condition (iii) holds.
Consider any $(x,q)\in R$, $l_{2}\in L_{2}$ and the transition $q\rTo^{l_{2}}_{2} p$ in $T_{\tau,\eta,\mu}(\Sigma)$. By definition of
$T_{\tau,\eta,\mu}(\Sigma)$
\begin{equation}
\left\Vert z-p\right\Vert \leq\eta/2, \label{b3}%
\end{equation}
where $z=\mathbf{x}(\tau,q,l_{2})\in Q_{1}$.
Choose $l_{1}=l_{2}\in L_{1}$ and consider now the transition \mbox{$x\rTo^{l_{1}}_{1} y$}
 in $T_{\mathcal{U}_{\tau}}(\Sigma)$. Since
$\Sigma$ is $\delta$--ISS and by (\ref{b3}) and (\ref{cond2}), the 
chain of inequalities in (\ref{gio21}) holds.
Thus $(y,p)\in R$, which completes the proof.
\end{proof}
%
%
%
\section{Symbolic control design for a pendulum\label{sec8}}
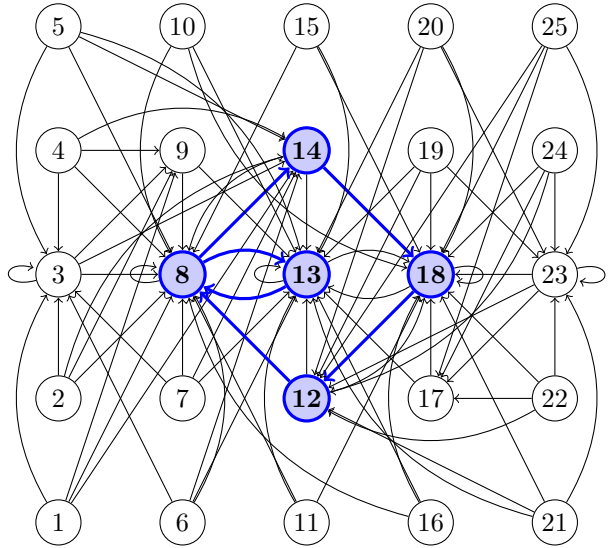
\begin{figure}[h]
\begin{center}
\begin{tikzpicture}[->,
shorten >=0.1pt,%
auto,node distance=1.65cm,thick,
inner sep=0.1 pt ,bend angle=30]
\tikzstyle{every state}=[minimum size=6mm]
\node[state, thin] (q_1) {1};
\node[state, thin] (q_2) [ above of=q_1] {2};
\node[state, thin] (q_3) [ above of=q_2] {3};
\node[state, thin] (q_4) [ above of=q_3] {4};
\node[state, thin] (q_5) [ above of=q_4] {5};
\node[state, thin] (q_6) [ right of=q_1] {6};
\node[state, thin] (q_7) [ above of=q_6] {7};
\node[state, draw=blue,very thick,fill=blue!20] (q_8) [ above of=q_7] {\bf{8}};
\node[state, thin] (q_9) [ above of=q_8] {9};
\node[state, thin] (q_10) [ above of=q_9] {10};
\node[state, thin] (q_11) [ right of=q_6] {11};
\node[state, draw=blue,very thick,fill=blue!20] (q_12) [ above of=q_11] {\bf{12}};
\node[state, draw=blue,very thick,fill=blue!20] (q_13) [ above of=q_12] {\bf{13}};
\node[state, draw=blue,very thick,fill=blue!20] (q_14) [ above of=q_13] {\bf{14}};
\node[state, thin] (q_15) [ above of=q_14] {15};
\node[state, thin] (q_16) [ right of=q_11] {16};
\node[state, thin] (q_17) [ above of=q_16] {17};
\node[state, draw=blue,very thick,fill=blue!20] (q_18) [ above of=q_17] {\bf{18}};
\node[state, thin] (q_19) [ above of=q_18] {19};
\node[state, thin] (q_20) [ above of=q_19] {20};
\node[state, thin] (q_21) [ right of=q_16] {21};
\node[state, thin] (q_22) [ above of=q_21] {22};
\node[state, thin] (q_23) [ above of=q_22] {23};
\node[state, thin] (q_24) [ above of=q_23] {24};
\node[state, thin] (q_25) [ above of=q_24] {25};
\path[->, thin] 
(q_1) 
edge [bend left] node {} (q_3)
edge node {} (q_8)
edge node {} (q_9)
edge node {} (q_14)
(q_2) 
edge node {} (q_3)
edge node {} (q_8)
edge node {} (q_9)
edge [bend left] node {} (q_14)
(q_3) 
edge [loop left] node {} (q_3)
edge node {} (q_8)
edge node {} (q_9)
edge node {} (q_14)
(q_4) 
edge node {} (q_3)
edge node {} (q_8)
edge node {} (q_9)
edge [bend left] node {} (q_14)
(q_5) 
edge [bend right] node {} (q_3)
edge node {} (q_8)
edge node {} (q_14)
edge [bend left] node {} (q_13)
(q_6) 
edge node {} (q_3)
edge [bend right] node {} (q_8)
edge node {} (q_13)
edge node {} (q_14)
(q_7) 
edge node {} (q_3)
edge node {} (q_8)
edge node {} (q_13)
edge node {} (q_14)
(q_8) 
edge [loop left] node {} (q_8)
edge [bend left, blue, very thick] node {} (q_13)
edge [blue, very thick] node {} (q_14)
(q_9) 
edge node {} (q_8)
edge node {} (q_13)
(q_10) 
edge [bend right] node {} (q_8)
edge node {} (q_13)
edge [bend right] node {} (q_18)
(q_11) 
edge node {} (q_8)
edge [bend left] node {} (q_13)
edge node {} (q_18)
(q_12) 
edge [blue, very thick] node {} (q_8)
edge node {} (q_13)
edge [bend right] node {} (q_18)
(q_13) 
edge [bend left, blue, very thick] node {} (q_8)
edge [loop left] node {} (q_13)
edge [bend left] node {} (q_18)
(q_14) 
edge [bend right] node {} (q_8)
edge node {} (q_13)
edge [blue, very thick] node {} (q_18)
(q_15) 
edge node {} (q_8)
edge [bend left] node {} (q_13)
edge node {} (q_18)
(q_16) 
edge [bend left] node {} (q_8)
edge node {} (q_13)
edge [bend left] node {} (q_18)
(q_17) 
edge node {} (q_13)
edge node {} (q_18)
(q_18) 
edge [blue, very thick] node {} (q_12)
edge [bend left] node {} (q_13)
edge [loop right] node {} (q_18)
(q_19) 
edge node {} (q_12)
edge node {} (q_13)
edge node {} (q_18)
edge node {} (q_23)
(q_20) 
edge [bend left] node {} (q_18)
edge node {} (q_12)
edge node {} (q_13)
edge node {} (q_23)
(q_21) 
edge node {} (q_12)
edge [bend left] node {} (q_13)
edge node {} (q_18)
edge [bend right] node {} (q_23)
(q_22) 
edge [bend left] node {} (q_12)
edge node {} (q_17)
edge node {} (q_18)
edge node {} (q_23)
(q_23) 
edge node {} (q_12)
edge node {} (q_18)
edge node {} (q_17)
edge [loop right] node {} (q_23)
(q_24) 
edge [bend left] node {} (q_12)
edge node {} (q_18)
edge node {} (q_17)
edge node {} (q_23)
(q_25) 
edge node {} (q_12)
edge node {} (q_17)
edge node {} (q_18)
edge [bend left] node {} (q_23);
\end{tikzpicture}
\end{center}
\caption{Symbolic model $T_{2,0.4,1.5\cdot 10^{-4}}(\Sigma)$ associated with the control system $\Sigma$ of (\ref{example}). 
\label{figure1}%
A state $(\eta\,i,\eta\,j)$ in $T_{2,0.4,1.5\cdot 10^{-4}}(\Sigma)$ with $i,j=-2,-1,0,1,2$ corresponds to the state $5\,(i+2)+j+3$ in the above picture.}
\end{figure}

One of the simplest mechanical control systems studied in the literature is the pendulum which can be described by:
\begin{equation}
\Sigma:\left\{
\begin{array}{clrr}
& \dot{x}_{1}=x_{2},\\
& \dot{x}_{2}=-\frac{g}{l}\sin\,x_{1}-\frac{k}{m}x_{2}+u,\\
\end{array}
\right.
\label{example}
\end{equation}
where $x_{1}$ and $x_{2}$ are the angular position and velocity of the point mass, $u$ is the torque which represents the control variable, $g=9.8$ is the gravity acceleration, $l=5$ is the length of the rod, $m=0.5$ is the mass and $k=3$ is the coefficient of friction. All constants and variables in system $\Sigma$ are expressed in the International System.
We assume that $u\in U=[-1.5,1.5]$ and that control inputs of $\Sigma$ are piecewise--constant. For simplicity we work on the subset $X=[-1,1]\times[-1,1]$ of the state space of $\Sigma$.\\
In order to apply Theorem \ref{Th:main5} we need to check if system $\Sigma$ is $\delta$--ISS. 
Consider the function $V:\mathbb{R}^{2}\times\mathbb{R}^{2}\rightarrow\mathbb{R}_{0}^{+}$ defined by:
\[
V(x,y)=\frac{1}{2}(x-y)'\left[
\begin{array}
{clrr}
\frac{1}{2}\left(\frac{k}{m}\right)^2 & \frac{1}{2}\frac{k}{m}\\
\frac{1}{2}\frac{k}{m} & \frac{1}{2}
\end{array}
\right](x-y).
\]
It is possible to show that $V$ satisfies condition (i) of Definition \ref{LyapunovGAS} with $\alpha_{1}(r)=0.49\,r^2$ and \mbox{$\alpha_{2}(r)=18.51\,r^2$}. 
Moreover, by defining for any $z_{1},z_{2}\in \mathbb{R}$,
\begin{equation}
\zeta(z_{1},z_{2})=(\sin(z_{1})-\sin(z_{2}))/(z_{1}-z_{2}),\notag\\
\end{equation}
one obtains $\zeta_{min}=\min_{z_{1},z_{2}\in [-1,1]}\zeta(z_{1},z_{2})=0.84$ and $\zeta_{max}=\max_{z_{1},z_{2}\in [-1,1]}\zeta(z_{1},z_{2})=1$ and hence:
%
\begin{eqnarray}
& &
\frac{\partial V}{\partial x}f(x,u)+\frac{\partial V}{\partial y}f(y,v) = 
-\frac{1}{2}\frac{k}{m}\frac{g}{l}\zeta(x_{1},y_{1})(x_{1}-y_{1})^{2}\notag\\
& &
-\frac{g}{l}\zeta(x_{1},y_{1})(x_{1}-y_{1})(x_{2}-y_{2})-\frac{1}{2}\frac{k}{m}(x_{2}-y_{2})^2\notag\\
& &
+\left(\frac{1}{2}\frac{k}{m}(x_{1}-y_{1})+x_{2}-y_{2}\right)(u-v)\notag\\
& &\leq -\frac{1}{2}a\Vert x-y \Vert_{2}^{2}+b|u-v|,
\label{Volta}
\end{eqnarray}
where $a=\frac{k}{m} \min\left\{\frac{g}{l}\,\zeta_{min},1\right\}-\frac{g}{l}\zeta_{max}=4.04>0$, \mbox{$b=(2+\frac{k}{m})=8>0$}.
Hence, condition (iii) of Definition \ref{LyapunovGAS} is satisfied with $\rho(r)=a\,r^2$ and $\sigma(r)=b\,r$, and $V$ is a $\delta$--ISS Lyapunov function for $\Sigma$. By Theorem \ref{TH-IISS} we conclude that the control system $\Sigma$ is $\delta$--ISS. 
Using inequality (\ref{Volta}), the definition of $V$ and the comparison lemma \cite{Khalil}, it is possible to show that for any $x,y\in X$, any $u,v\in \mathcal{U}$ and any time $t\in\mathbb{R}^{+}_{0}$:
\begin{equation}
\left\Vert \mathbf{x}(t,x,u)-\mathbf{x}(t,y,v)\right\Vert
\leq 
\beta(\left\Vert x-y \right\Vert,t)+\gamma(\Vert u-v \Vert_{\infty}),\notag
\end{equation}
where $\beta(r,s):=6.17\,e^{-2.08\,s}r$ and $\gamma(r):=\sqrt{3.96\,r}$ for any \mbox{$r,s\in \mathbb{R}$}. Functions $\beta$ and $\gamma$ are respectively $\mathcal{KL}$ and $\mathcal{K}_{\infty}$ functions and thus inequality (\ref{deltaISS}) is satisfied.
We now have all the ingredients to apply Theorem \ref{Th:main5}. Condition (\ref{cond2}) becomes:
\begin{equation}
6.17\,e^{-2.08\,\tau}\varepsilon+\sqrt{3.96\,\mu}+\eta/2\leq\varepsilon.
\label{ineq}%
\end{equation}
For a precision $\varepsilon=0.25$ we can choose $\eta=0.4$, \mbox{$\tau=2$} and $\mu=1.5 \cdot  10^{-4}$ so that inequality (\ref{ineq}) is satisfied. 
The resulting transition system:
\begin{equation}
T_{\tau,\eta,\mu}(\Sigma)=(Q_{2},L_{2},\rTo_{2},O_{2},H_{2}),
\label{Texample}
\end{equation}
 is defined by:
\begin{itemize}
\item $Q_{2}=\{-2\eta,-\eta,0,\eta,2\eta\}\times\{-2\eta,-\eta,0,\eta,\,2\eta\}$;
\item $L_{2}=[U]_{1.5 \cdot  10^{-4}}$;
\item $\rTo_{2}$ is depicted in Figure \ref{figure1};
\item $O_{2}=X$;
\item $H_{2}=\imath : Q_{2} \hookrightarrow O_{2}$,
\end{itemize}
and shown in Figure \ref{figure1} where the transition relation $\rTo_{2}$ has been obtained by numerically integrating the trajectories of $\Sigma$.\\

\begin{figure}[ptb]
\begin{center}
\includegraphics[scale=0.35]{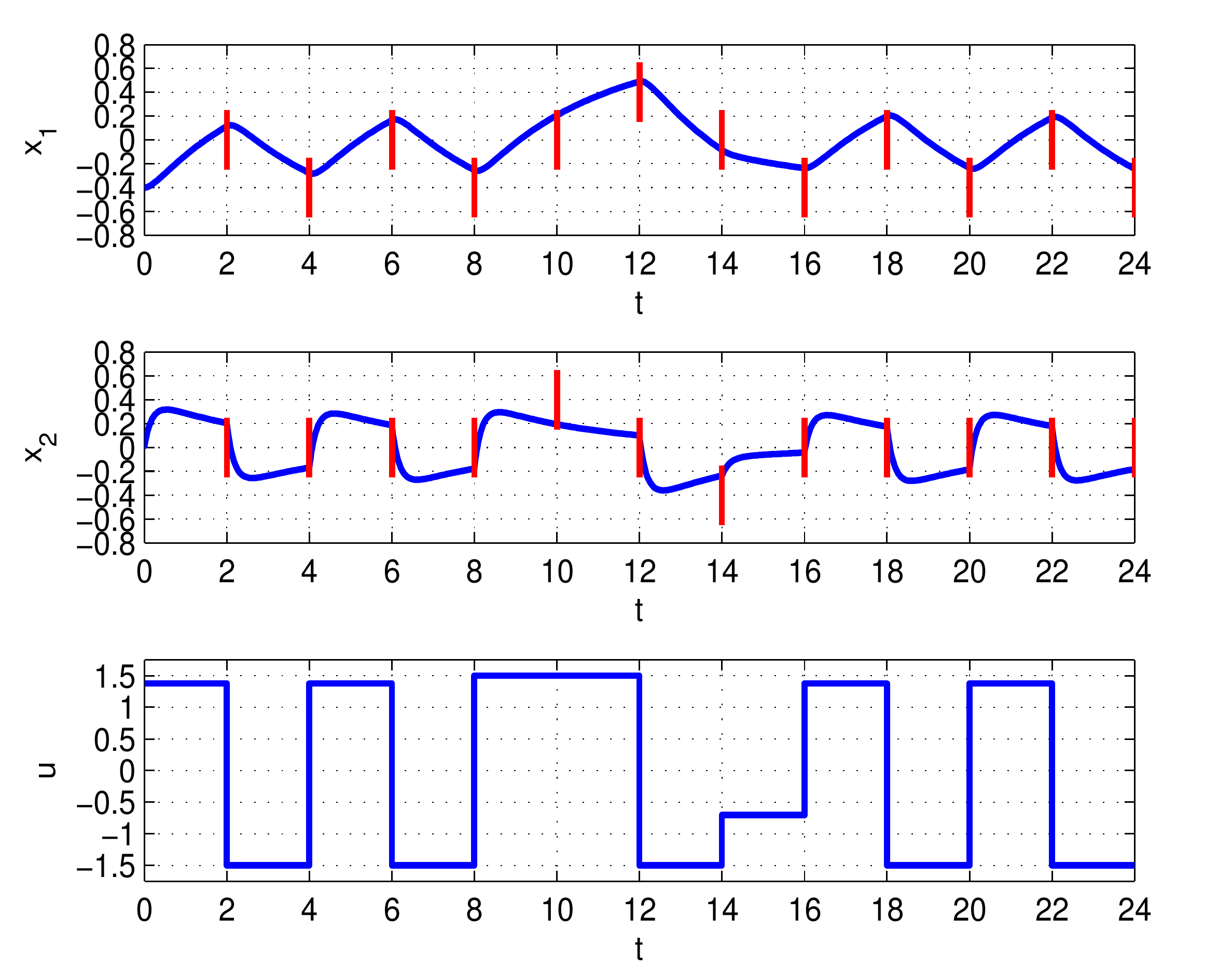}
\end{center}
\caption{Upper and medium panels: trajectory of $(x_{1},x_{2})$, with initial condition $(-\eta,0)$ and control strategy synthesized on $T_{2,0.4,1.5\cdot 10^{-4}}(\Sigma)$. Vertical intervals marked in red represent the precision $\varepsilon=0.25$ that we require.  Lower panel: Control strategy synthesized on $T_{2,0.4,1.5\cdot 10^{-4}}(\Sigma)$.}%
\label{fig3}%
\end{figure}

We now illustrate the use of the symbolic model (\ref{Texample}) for controller synthesis. Suppose that our objective is to design a controller enforcing an alternation between two different periodic motions denoted by $P_1$ and $P_2$. Periodic motion $P_1$ requires  the state of $\Sigma$ to cycle between $(-\eta,0)$ and $(0,0)$ while periodic motion $P_{2}$ requires the state to cycle between $(-\eta,0)$ and $(\eta,0)$. The control objective is then the design of a controller that enforces system $\Sigma$ to satisfy a specification $P$ requiring the execution of the sequence of  periodic motions $P_{1},P_{1},P_{2},P_{1},P_{1}$. 
This specification is a simple illustration of more complex control objectives that typically require different sequencing of actions in response to exogenous events such as faults or to events triggered by the violation of certain thresholds on the continuous state. This kind of specifications will naturally result in a hybrid controller combining the continuous inputs necessary to drive the continuous state with the discrete logic responsible for executing the right sequence of actions in response to different conditions.
A control strategy for periodic motions $P_{1}$ and $P_2$ can be obtained by performing a simple search\footnote{States 
$(-\eta,0)$, $(0,0)$, $(0,\eta)$, $(\eta,0)$ and $(0,-\eta)$ involved in the specifications $P_{1}$ and $P_{2}$, correspond respectively to states $8$, $13$, $14$, $18$ and $12$ in Figure \ref{figure1}.} on $T_{2,0.4,1.5\cdot 10^{-4}}(\Sigma)$ or by using standard methods in the context of supervisory control \cite{RamWonI} or algorithmic approaches to game theory \cite{AVW03}. One possible solution enforcing $P_1$ is:
\[
(-\eta,0)\rTo^{1.38}(0,0)\rTo^{-1.5}(-\eta,0),
\]
and for $P_2$ is:
\[
(-\eta,0)\rTo^{1.5}(0,\eta)\rTo^{1.5}(\eta,0)\rTo^{-1.5}(0,-\eta)\rTo^{-0.71}(-\eta,0).
\]
A control strategy that enforces the specification $P$ can be obtained by concatenating the trajectories associated with $P_{1}$, $P_{1}$, $P_{2}$, $P_{1}$ and $P_{1}$, resulting in:
\begin{eqnarray}
& &(-\eta,0)\rTo^{1.38}(0,0)\rTo^{-1.5}(-\eta,0)\rTo^{1.38}(0,0)\rTo^{-1.5}(-\eta,0)\notag\\
& & \rTo^{1.5}(0,\eta)\rTo^{1.5}(\eta,0)\rTo^{-1.5}(0,-\eta)\rTo^{-0.71}(-\eta,0)\notag\\
& & \rTo^{1.38}(0,0)\rTo^{-1.5}(-\eta,0)\rTo^{1.38}(0,0)\rTo^{-1.5}(-\eta,0).\notag
\end{eqnarray}
Since by Theorem \ref{Th:main5}, $T_{2,0.4,1.5\cdot 10^{-4}}(\Sigma)$ is \mbox{$0.25$--bisimilar} to $T_{\mathcal{U}_{2}}(\Sigma)$, the notion of approximate bisimulation guarantees that the controller synthesized on $T_{2,0.4,1.5\cdot 10^{-4}}(\Sigma)$, will enforce the desired behavior on $\Sigma$ with an error of at most $0.25$. Figure \ref{fig3} shows the evolution of the state variables of $\Sigma$, when applying such control strategy. It is easy to see that at each time $i\,\tau$ with $i=1,...,12$ the state variables $x_{1}$ and $x_{2}$ are within the interval marked in red, which represents the desired precision $\varepsilon=0.25$. For example, at time $t=2\,\tau=4$ the angular position $x_{1}$ of system $\Sigma$ is in the interval $-\eta+[-\varepsilon,\varepsilon]=[-0.65,-0.15]$, as required by $P_{1}$ and the approximation error $\varepsilon$. 
Although we could have designed continuous controllers enforcing $P_1$ and $P_2$ and then devise a switching logic enforcing specification $P$, as is currently done in practice, we could not guarantee what would happen to the closed loop system due to the difficulty in analyzing the combination of continuous controllers with switching logic (see e.g. \cite{liberzon2003}). 
On the contrary, the methodology that we propose offers a systematic controller design process that requires reduced user intervents.
\section{Discussion\label{sec9}}
The work presented in this paper compares as follows with the available results of the research lines recalled in the introduction.

\textbf{Simulation/bisimulation:}
The results in this paper follow the research line of \cite{Paulo_HSCC07} and provide important generalizations:\\
(i) The definition of the symbolic model in \cite{Paulo_HSCC07} relies on an (arbitrary) a--priori choice of control inputs, while the symbolic model in (\ref{T1}) captures the effect of \textit{any} measurable control input;\\
(ii) The approximation notion employed in \cite{Paulo_HSCC07} is approximate simulation\footnote{
We recall from \cite{AB-TAC07} that an \mbox{$\varepsilon$--approximate} simulation relation from $T_{1}$ to $T_{2}$ is a relation $R$ which satisfies conditions (i) and (ii) in Definition \ref{ASR}.} while the results in this paper guarantee the stronger notion of approximate \emph{bisimulation}.\\
These generalizations are quite important from the controller synthesis point of view. The main drawback of the results in \cite{Paulo_HSCC07} is that if a controller fails to exist for the symbolic model, nothing can be concluded regarding the existence of a controller for the original control system. 
Our results guarantee, instead, that given a control system and a specification, a controller exists for the original model if and only if a controller exists for the symbolic model.
Notice that while $\delta$--GAS implies asymptotic stabilizability as employed in Theorem 2 of \cite{Paulo_HSCC07}, the converse is not true in general\footnote{In fact the converse is true in the case of linear control systems.}. Furthermore even if a feedback control law rendering the closed--loop system $\delta$--GAS were found, if the input space of the control system is bounded, there is no guarantee that such feedback would satisfy the input constraints.\\
The results in this paper share similar ideas with the ones in \cite{Girard_HSCC07} that considers discrete--time linear control systems. 
When we regard discrete--time control systems as the time discretization of continuous--time control systems, Theorem \ref{Th:main5} extends Theorem 4 of \cite{Girard_HSCC07} in two directions:\\
(i) by enlarging the class of control systems from linear to nonlinear;\\
(ii) by enlarging the class of input signals from piecewise--constant to measurable.\\
When specializing results of this paper to the class of \textit{linear} control systems, conditions of Theorems \ref{Th_main} and \ref{Th:main5} simplify. In fact given a linear control system:
\[
\dot{x}=Ax+Bu,\text{      }x\in \mathbb{R}^{n},\text{      }u\in U\subseteq\mathbb{R}^{m},
\]
the notions of $\delta$--GAS and $\delta$--ISS reduce to asymptotic stability of matrix $A$ and functions $\beta$ and $\gamma$ appearing in inequalities (\ref{deltaUGAS}) and (\ref{deltaISS}) can be chosen as:
\begin{eqnarray}
\beta(r,s)=\Vert e^{As} \Vert r; 
&\text{                      }& 
\gamma(r)=\left( \Vert B \Vert \int_{0}^{\infty}\Vert e^{As}\Vert ds\right) r,
\label{betagamma}
\end{eqnarray}
where $\Vert e^{As} \Vert$ denotes the infinity norm of the matrix\footnote{For $M=\{m_{ij}\}\in\mathbb{R}^{n\times m}$, $\Vert M\Vert:=\max_{1\leq i\leq m}
{\textstyle\sum_{j=1}^{n}}
|m_{ij}|$.} $e^{As}$.
The use of explicit expressions in (\ref{betagamma}) for $\beta$ and $\gamma$ simplifies indeed the search of parameters 
$\tau$, $\eta$ and $\mu$ satisfying conditions of Theorems \ref{Th_main} and \ref{Th:main5}, and hence the construction of symbolic models in (\ref{T1}) and (\ref{T2}). Furthermore, in contrast to the nonlinear case, the construction of the symbolic models can be performed even for non-constant inputs. This can be done by using results on polytopic approximation of reachable sets for linear control systems 
(see e.g. \cite{Varaiya:98}, \cite{Girard_HSCC05}) with compact input space. 
It is known from \cite{Varaiya:98} that for any desired precision \mbox{$\nu\in \mathbb{R}^{+}$}, the reachable set $\mathcal{R}(\tau,q)$ of (\ref{reach}) can be approximated by a polytope $P(\tau,q)$, so that
$\mathbf{d}_{h}(P(\tau,q),\mathcal{R}(\tau,q))\leq \nu$,
where $\mathbf{d}_{h}$ is the Hausdorff pseudo--metric\footnote{We recall that for any $X_{1},X_{2}\subseteq \mathbb{R}^{n}$,  \mbox{$\mathbf{d}_{h}(X_{1},X_{2}):=\max\{\vec{\mathbf{d}}_{h}(X_{1},X_{2}),\vec{\mathbf{d}}_{h}(X_{2},X_{1})\}$},
where $\vec{\mathbf{d}}_{h}(X_{1},X_{2}):=\sup\nolimits_{x_{1}\in X_{1}}\inf
\nolimits_{x_{2}\in X_{2}}\mathbf{d}(x_{1},x_{2})$.} induced by the metric $\mathbf{d}$. The countable set $\mathcal{P}_{\mu}(\tau,q)$, 
can then be reformulated in terms of $P(\tau,q)$ rather than of $\mathcal{R}(\tau,q)$, as follows:
\[
\mathcal{P}_{\mu}(\tau,q):=\{y\in\lbrack\mathbb{R}^{n}]_{\mu}:\exists
z\in P(\tau,q)\text{ s.t. }\Vert y-z\Vert\leq\mu/2\}\text{.}%
\]
The symbolic model in (\ref{T1}) can be adapted to the case of linear systems by defining the set $L_{2}(q)$ by:
\begin{equation}
L_{2}(q):=\mathcal{P}_{\mu}(\tau,q)
\label{L2q}
\end{equation}
and the transition relation $\rTo_{2}$ by: 
\begin{eqnarray}
q \rTo^{l}_{2} p,
& 
\text{ if     }
& 
\left\Vert p-\mathbf{x}(\tau,q,0)-l\right\Vert \leq\eta/2.
\label{Trelation}
\end{eqnarray}
Since the sets $P(\tau,q)$ and $L_{2}(q)$ can be computed the symbolic model (\ref{T1}) with $L_{2}(q)$ given by (\ref{L2q}) and $\rTo_{2}$ given by (\ref{Trelation}), can be constructed.
Finally condition (\ref{cond}) of Theorem \ref{Th_main} can be adapted to this case, resulting in 
$\Vert e^{A\tau} \Vert \varepsilon+\nu+\mu+\eta/2\leq\varepsilon$.

\textbf{Quantized control systems:}
In \cite{BMP02,BMP06} finite abstractions of quantized control systems are studied. In particular, conditions on the systems parameters and on the input set are found so that the resulting abstraction is characterized by a lattice structure in the set $\mathcal{R}$ of reachable states. Our results ensure, under the $\delta$--ISS assumption, existence of a lattice approximating $\mathcal{R}$, independently from the system parameters and input set. 
More precisely a direct consequence of Theorem \ref{Th:main5} is that if a digital control system $\Sigma$ is $\delta$--ISS then any state $x\in\mathcal{R}$ can be approximated with any desired precision $\varepsilon\in\mathbb{R}^{+}$, by a (symbolic) state $q\in [\mathbb{R}^{n}]_{\varepsilon}$ so that $\Vert x-q \Vert \leq \varepsilon/2$. However, while our results guarantee to \emph{approximate} $\mathcal{R}$ by the lattice $[\mathbb{R}^{n}]_{\varepsilon}$ with any (arbitrarily small) precision $\varepsilon\in\mathbb{R}^{+}$, results established in \cite{BMP02,BMP06} guarantee that $\mathcal{R}$ \emph{is exactly} a lattice.

\textbf{Qualitative reasoning and Stochastic automata:}
Symbolic models have been also proposed in the framework of qualitative reasoning (see
e.g. \cite{QR_HSCC03,QRBook}) and in the stochastic automata based abstraction of \cite{Lunze2001,LunzeBook}. 
In both approaches the proposed models are characterized by a ``completeness'' property under which, any trajectory of the control system can be mimicked by a trajectory of the proposed symbolic models. On the other hand, for any trajectory of the symbolic models there may not exist a corresponding matching trajectory in the control systems. In both approaches no stability assumptions are needed to ensure the completeness property. An interpretation in terms of bisimulation theory, is that these results guarantee existence of a surjective exact simulation relation\footnote{An exact simulation relation is an $\varepsilon$--approximate simulation relation with $\varepsilon=0$.} from the control systems to the symbolic models. 
However, analogously to the results in \cite{Paulo_HSCC07} the main drawback of these approaches is that if a controller fails to exist for the proposed symbolic models, nothing can be concluded regarding the existence of a controller for the original control system. 
As pointed out before, this drawback can be overcome by considering a notion of approximate bisimulation, whose existence is ensured by $\delta$--ISS of the control system (see Theorem \ref{Th:main5}).

The results in Section \ref{sec6} provide a first step towards the effective computation of symbolic models for digital control systems. However, further work is required towards the design of efficient algorithms for constructing the symbolic model proposed in (\ref{T2}). In particular, the main critical issues 
are related with:\\
(i) the choice of parameters $\tau,\eta,\mu$, which translates, by inequality (\ref{cond2}), in finding 
 a $\delta$--ISS Lyapunov function for the control system;\\
(ii) the cardinality of $Q_{2}$ and $L_{2}$, which increases exponentially with the dimension of the state and input spaces of the control system.\\
The computation of $\delta$--ISS Lyapunov functions is in general a hard task. However,
one can resort to numerical tools available in the literature, as for example the one proposed in \cite{SOS1}. 
Furthermore, a way for mitigating the exponential grow in the sizes of $Q_{2}$ and $L_{2}$ is to adapt techniques 
from on–-the-–fly verification of transition systems \cite{onthefly2} to the construction of the proposed symbolic models. This will be the object of future investigations.

\bibliographystyle{alpha}
\bibliography{biblio}

\end{document}